\documentclass[12pt]{article}
\usepackage{amsfonts, amssymb}
\begin{document}
\baselineskip 18pt \hoffset=-1cm \voffset=-1in
\renewcommand{\theequation}{\arabic{section}.\arabic{equation}}
\newtheorem{theorem}{Theorem}[section]
\newtheorem{lemma}{Lemma}[section]
\newtheorem{proposition}{Proposition}[section]
\newtheorem{corollary}{Corollary}[section]
\newtheorem{definition}{Definition}[section]
\newtheorem{algorithm}{Algorithm}[section]
\newtheorem{remark}{Remark}[section]
\newtheorem{exercise}{Exercise}[section]
\def\Proof{\noindent{\bf Proof:~}}
\def\blackslug{\hbox{\hskip 1pt \vrule width 4pt height 8pt
    depth 1.5pt \hskip 1pt}}
\def\QED{\quad\blackslug\lower 6.5pt\null\par}
\def\R{\mathbb{R}}
\def\X{\mathbb{R}^J}
\def\C{\mathbb{C}}
\def\H{{\cal H}}

\title {Alternating Minimization, Proximal Minimization and Optimization Transfer Are Equivalent }

\author{Charles L. Byrne and Jong Soo Lee \\ Department of
Mathematical Sciences \\ University of Massachusetts Lowell \\ Lowell,
MA 01854}
\date{December 9, 2015}
\maketitle

\begin{abstract}Let $X$ be an arbitrary nonempty set and $f:X\rightarrow \R$. The objective is to minimize $f(x)$ over $x\in X$. The iterative algorithms considered here are \lq\lq descent\rq\rq algorithms, so that $\{f(x^k)\}\downarrow \beta^*\geq -\infty$. We want $\beta^*=\beta\doteq \inf_{x\in X}f(x)$. 

In proximal minimization algorithms (PMA) we minimize
$f(x)+d(x,x^{k-1})$ to get $x^k$. The $d:X\times X\rightarrow \R_+$ is a \lq\lq distance\rq\rq  function, with $d(x,x)=0$, for all $x$. 

In majorization minimization (MM), also called optimization transfer, a second \lq\lq majorizing\rq\rq\, function $g(x|z)$ is introduced, with the properties $g(x|z)\geq f(x)$, for all $x$ and $z$ in $X$, and $g(x|x)=f(x)$. We then minimize $g(x|x^{k-1})$ to get $x^k$. 

Let $\Phi :X\times Y\rightarrow \R_+$, where $X$ and $Y$ are arbitrary nonempty sets. The objective in alternating minimization (AM) is to find ${\hat x}\in X$ and ${\hat y}\in Y$ such that 
$\Phi ({\hat x},{\hat y})\leq \Phi (x,y)$ for all $x\in X$ and $y\in Y$. For each $k$ we minimize $\Phi (x,y^{k-1})$ to get $x^{k-1}$ and then minimize $\Phi (x^{k-1},y)$ to get $y^k$. For each $x\in X$, let $y(x)\in Y$ be such that $\Phi (x,y)\geq \Phi (x,y(x))$, for all $y\in Y$; then $y^k=y(x^{k-1})$. Minimizing $\Phi (x,y)$ over all $x\in X$ and $y\in Y$ is equivalent to minimizing $f(x)\doteq \Phi (x,y(x))$ over all $x\in X$. With $d(x,z)\doteq\Phi (x,y(z))-\Phi (x,y(x))$, minimizing $\Phi (x,y^k)$ is equivalent to minimizing $f(x)+d(x,x^{k-1})$. Therefore, AM, MM, and PMA are equivalent. Each type of algorithm leads to a decreasing sequence $\{f(x^k)\}$.

New conditions on PMA that imply $\beta^*=\beta$ are given, which lead to new conditions on AM for the sequence $\{\Phi (x^k,y^k)\}$ to converge to $\inf_{x,y}\Phi (x,y)$. These conditions can then be translated into the language of MM. Examples are given of each type of algorithm and some open questions are posed.



\end{abstract}

\section{Introduction}\setcounter{equation}{0}
\setcounter{theorem}{0} \setcounter{lemma}{0}
\setcounter{proposition}{0} \setcounter{corollary}{0}
\setcounter{definition}{0} \setcounter{algorithm}{0}Let $X$ be an arbitrary nonempty set and $f:X\rightarrow \R$. The objective is to minimize $f(x)$ over $x\in X$. The iterative algorithms considered here are \lq\lq descent\rq\rq algorithms, so that $\{f(x^k)\}\downarrow \beta^*\geq -\infty$. We want $\beta^*=\beta\doteq \inf_{x\in X}f(x)$.

In proximal minimization algorithms (PMA) we minimize
$f(x)+d(x,x^{k-1})$ to get $x^k$. The $d:X\times X\rightarrow \R_+$ is a \lq\lq distance\rq\rq  function, with $d(x,x)=0$, for all $x$. In majorization minimization (MM), also called optimization transfer, a second \lq\lq majorizing\rq\rq\, function $g(x|z)$ is introduced, with the properties $g(x|z)\geq f(x)$, for all $x$ and $z$ in $X$, and $g(x|x)=f(x)$. We then minimize $g(x|x^{k-1})$ to get $x^k$. With
$$d(x,z)\doteq g(x|z)-f(x),$$ it is clear that MM is equivalent to PMA; alternating minimization (AM) algorithms appear to be more general.

Let $\Phi :X\times Y\rightarrow \R_+$, where $X$ and $Y$ are arbitrary nonempty sets. The objective in AM is to find ${\hat x}\in X$ and ${\hat y}\in Y$ such that 
$$\Phi ({\hat x},{\hat y})\leq \Phi (x,y),$$ for all $x\in X$ and $y\in Y$. For each $k$ we minimize $\Phi (x,y^{k-1})$ to get $x^{k-1}$ and then minimize $\Phi (x^{k-1},y)$ to get $y^k$. We have the following proposition:
\begin{proposition} The AM, PMA, and MM methods are equivalent. \label{equiv}\end{proposition}
\Proof We reformulate AM as a method for minimizing a function $f(x)$ of the single variable $x\in X$. For each $x\in X$, let $y(x)\in Y$ be such that $\Phi (x,y)\geq \Phi (x,y(x))$, for all $y\in Y$. Then minimizing $\Phi (x,y)$ over all $x\in X$ and $y\in Y$ is equivalent to minimizing $f(x)\doteq \Phi (x,y(x))$ over all $x\in X$.  
Every MM algorithm, and therefore every PMA, can be viewed as an application of alternating minimization: define $\Phi (x,z)\doteq g(x|z)$. Minimizing $g(x|x^{k-1})$ to get $x^k$ is equivalent to minimizing $\Phi (x,x^{k-1})$, while minimizing $g(x^k|z)$ is equivalent to minimizing $\Phi (x^k,z)$ and yields $z=x^k$.
\hfill \QED 
Note that $\Phi (x^{k-1},y^k)=f(x^{k-1})$. Then the sequence $\{f(x^k)\}$ is decreasing to some $\beta^*$.

Each of the algorithms we consider can be reformulated as minimizing some objective function $f(x)$ and can be described by saying that at each step we minimize
$$G_k(x)=f(x)+g_k(x),$$ where $g_k(x)\geq 0$ and $g_k(x^{k-1})=0$. Such methods are called {\it auxiliary-function} (AF) algorithms \cite{Byr08}. For AF algorithms we know that the sequence $\{f(x^k)\}$ is decreasing to some number $\beta^*\geq -\infty$. If an AF algorithm is in the subclass of SUMMA2 algorithms, then we know that $\beta^*=\beta \doteq\inf_x f(x)$. The Euclidean and Kullback-Leibler distances yield algorithms in the SUMMA2 class, and we suspect that the methods based on the Hellinger and Pearson $\phi^2$ distances are also in the SUMMA2 class. Conditions are presented that are sufficient for PMA to be in the SUMMA2 class, and therefore, for $\beta^*=\beta$ for AM, PMA, and MM algorithms. We also consider the use of alternating minimization of distances to obtain approximate solutions of systems of linear equations. The distances considered include the Euclidean, the Kullback-Leibler,  the Hellinger, and the Pearson $\phi^2$ distances. 

\section{Auxiliary-Function Methods in Optimization}\setcounter{equation}{0}
\setcounter{theorem}{0} \setcounter{lemma}{0}
\setcounter{proposition}{0} \setcounter{corollary}{0}
\setcounter{definition}{0} \setcounter{algorithm}{0}

Let $f:X\rightarrow \R$, where $X$ is an arbitrary nonempty set. In applications the set $X$ will have additional structure, but not always that of a Euclidean space; for that reason, it is convenient to impose no structure at the outset. An iterative procedure for minimizing $f(x)$ over $x\in X$ is called an {\it auxiliary-function} (AF) algorithm \cite{Byr08} if, at each step, we minimize
\begin{eqnarray}G_k(x)=f(x)+g_k(x),\end{eqnarray} where $g_k(x)\geq 0$, and $g_k(x^{k-1})=0$. It follows easily that the sequence $\{f(x^k)\}$ is decreasing, so $\{f(x^k)\}\downarrow \beta^*\geq -\infty $.
We want more, however; we want $\beta^*=\beta \doteq \inf_{x\in X}f(x)$. To have this we need to impose an additional condition on the auxiliary functions $g_k(x)$; the SUMMA Inequality \cite{Byr08} is one such additional condition. 

\subsection{The SUMMA Class} We say that an AF algorithm is in the SUMMA class if the SUMMA Inequality holds for all $x$ in $X$:
\begin{eqnarray}G_k(x)-G_k(x^k)\geq g_{k+1}(x).\label{SUMMAinequality}\end{eqnarray}
One consequence of the SUMMA Inequality is
\begin{eqnarray} g_k(x)+f(x)\geq g_{k+1}(x)+ f(x^k),\label{SI1}\end{eqnarray} for all $x\in X$. It follows from this that $\beta^*=\beta$. If this were not the case, then there would be $z\in X$ with 
$$f(x^k)\geq \beta^*> f(z)$$ for all $k$. The sequence $\{g_k(z)\}$ would then be a decreasing sequence of nonnegative terms with the sequence of its successive differences bounded below by $\beta^*-f(z)>0$.

There are many iterative algorithms that satisfy the SUMMA Inequality \cite{Byr08}, such as barrier-function methods \cite{FM90}, and are therefore in the SUMMA class. However, some important methods that are not in this class still have $\beta^*=\beta$; one example is the proximal minimization method of Auslender and Teboulle \cite{AT06}.  This suggests that the SUMMA class, large as it is, is still unnecessarily restrictive. This leads us to the definition of the SUMMA2 class.

\subsection{The SUMMA2 Class} An iterative algorithm for minimizing $f:X\rightarrow \R$ is said to be in the SUMMA2 class if, for each sequence $\{x^k\}$ generated by the algorithm, there are functions $h_k:X\rightarrow \R_+$ such that, for all $x\in X$, we have
\begin{eqnarray}h_k(x)+f(x)\geq h_{k+1}(x)+f(x^k).\label{SI2}\end{eqnarray} Any algorithm in the SUMMA class is in the SUMMA2 class; use $h_k=g_k$. As in the SUMMA case, we must have $\beta^*=\beta$, since otherwise the successive differences of the sequence $\{h_k(z)\}$ would be bounded below by $\beta^*-f(z)>0$. It is helpful to note that the functions $h_k$ need not be the $g_k$, and we do not require that $h_k(x^{k-1})=0$. The proximal minimization method of Auslender and Teboulle is in the SUMMA2 class, as is the expectation maximization maximum likelihood (EMML) algorithm \cite{SV82, VSK85, Byr93}.

\section{PMA is MM}\setcounter{equation}{0}
\setcounter{theorem}{0} \setcounter{lemma}{0}
\setcounter{proposition}{0} \setcounter{corollary}{0}
\setcounter{definition}{0} \setcounter{algorithm}{0} In proximal minimization algorithms (PMA) we minimize
\begin{eqnarray}f(x)+d(x,x^{k-1})\end{eqnarray} to get $x^k$. Here $d(x,z)\geq 0$ and $d(x,x)=0$, so we say that $d(x,z)$ is a distance.

In \cite{CZL14} the authors review the use, in statistics, of \lq\lq majorization minimization\rq\rq\, (MM), also called \lq\lq optimization transfer\rq\rq . In numerous papers \cite{EF99, AFBH06} Jeff Fessler and his colleagues use the terminology \lq\lq surrogate-function minimization\rq\rq\, to describe optimization transfer. The objective is to minimize $f:X\rightarrow \R$. In MM methods a second \lq\lq majorizing\rq\rq\, function $g(x|z)$ is introduced, with the properties $g(x|z)\geq f(x)$, for all $x$ and $z$ in $X$, and $g(x|x)=f(x)$. We then minimize $g(x|x^{k-1})$ to get $x^k$. Defining
$$d(x,z)\doteq g(x|z)-f(x),$$ it is clear that $d(x,z)$ is a distance and so MM is equivalent to PMA.

\section{PMA with Bregman Distances (PMAB)} Let $f:\R^J\rightarrow \R$ and $h:\R^J\rightarrow \R$ both be convex and differentiable. Let 
$$D_h(x,z)\doteq h(x)-h(z)-\langle \nabla h(z),x-z\rangle $$ be the Bregman distance associated with $h$.  At the $k$th step of a proximal minimization algorithm with Bregman distance (PMAB) we minimize
\begin{eqnarray}G_k(x)=f(x)+D_h(x,x^{k-1})=f(x)+h(x)-h(x^{k-1})-\langle \nabla h(x^{k-1}),x-x^{k-1}\rangle\end{eqnarray} to get $x^k$. It was shown in \cite{Byr08} that 
$$G_k(x)-G_k(x^k)=D_f(x,x^k)+D_h(x,x^k)\geq D_h(x,x^k)=g_{k+1}(x),$$ so that all PMAB are in the SUMMA class.

In order to minimize $G_k(x)$ we need to solve the equation
\begin{eqnarray}0=\nabla f(x)+\nabla h(x)-\nabla h(x^{k-1})\end{eqnarray} for $x=x^k$; generally, this is not easy. Here is a \lq\lq trick\rq\rq\, that can be used to simplify the calculations. Select a function $g$ so that $h\doteq g-f$ is convex and differentiable and so that the equation
\begin{eqnarray}0=\nabla g(x)-\nabla g(x^{k-1})+\nabla f(x^{k-1})\end{eqnarray} is easily solved. As an example, we use this \lq\lq trick\rq\rq\, to derive a gradient descent algorithm and the Landweber algorithm.

\section{Gradient Descent and the Landweber Algorithm} Suppose that we want to minimize a convex differentiable function $f:\R^J\rightarrow \R$. If the gradient of $f$, $\nabla f$, is a $L$-Lipschitz continuous operator, that is, if
$$\|\nabla f(x)-\nabla f(z)\|\leq L \|x-z\|,$$ then the function
$$ h(x)\doteq g(x)-f(x)=\frac{1}{\gamma}\|x\|^2-f(x)$$ is convex, for $0<\gamma \leq 1/L$. For each $k$ we minimize
$$G_k(x)=f(x)+\frac{1}{\gamma}\|x-x^{k-1}\|^2-D_f(x,x^{k-1})$$ to get $x^k$. We then have
$$x^k=x^{k-1}-\gamma \nabla f(x^{k-1}),$$ which is a gradient descent algorithm. As a special case we get Landweber's algorithm.

Suppose we want to find a minimizer of the function $f(x)=\|Ax-b\|^2$, where $A$ is a real $I$ by $J$ matrix. Let $g(x)=\frac{1}{\gamma}\|x\|^2$, for some $\gamma$ in the interval $(0,\frac{1}{L})$, where $L=\rho (A^TA)$, the largest eigenvalue of the matrix $A^TA$. Then the function $h\doteq g-f$ is convex and differentiable. We have
\begin{eqnarray}D_f(x,y)=\|Ax-Ay\|^2,\end{eqnarray} so that
\begin{eqnarray}D_h(x,y)=\frac{1}{\gamma}\|x-y\|^2-\|Ax-Ay\|^2.\end{eqnarray} At the $k$th step we differentiate
\begin{eqnarray}\|Ax-b\|^2+\frac{1}{\gamma}\|x-x^{k-1}\|^2-\|Ax-Ax^{k-1}\|^2,\end{eqnarray} to obtain
\begin{eqnarray}0=A^T(Ax-b)+\frac{1}{\gamma}(x-x^{k-1})-A^T(Ax-Ax^{k-1}),\end{eqnarray} so that
\begin{eqnarray}x^k=x^{k-1}-\gamma A^T(Ax^{k-1}-b).\end{eqnarray} This is the iterative step of Landweber's algorithm. The sequence $\{x^k\}$ converges to a minimizer $x^*$ of $f(x)$, and $x^*$ minimizes $\|{\hat x}-x^0\|$ over all ${\hat x}$ that minimize $\|Ax-b\|$.

In \cite{Byr14} this same \lq\lq trick\rq\rq was used to obtain an elementary proof of convergence of the forward-backward-splitting algorithm \cite{CW05}.

\section{The Quadratic Upper Bound Principle} In \cite{BL88} the authors introduce the {\it quadratic upper bound principle} as a method for obtaining a majorizing function in optimization transfer. The objective is to minimize the function $f:\R^J\rightarrow \R$. If $f$ is twice continuously differentiable, then, for any $x$ and $z$, we have, according to the extended Mean Value Theorem,
\begin{eqnarray}f(x)=f(z)+\langle \nabla f(z),x-z\rangle +\frac{1}{2}(x-z)^T\nabla ^2f(w)(x-z),\end{eqnarray} for some $w$ on the line segment connecting $x$ and $z$. If there is a positive-definite matrix $B$ such that $B-\nabla^2f(w)$ is positive-definite for all $w$, then we have
\begin{eqnarray}f(x)\leq f(z)+\langle \nabla f(z),x-z\rangle +\frac{1}{2}(x-z)^TB(x-z).\end{eqnarray}
Then we have $g(x|z)\geq f(x)$, for all $x$ and $z$, where
\begin{eqnarray}g(x|z)\doteq f(z)+\langle \nabla f(z),x-z\rangle +\frac{1}{2}(x-z)^TB(x-z).\end{eqnarray} The iterative step is now to minimize $g(x|x^{k-1})$ to get $x^k$.

The iterative step is equivalent to minimizing
\begin{eqnarray}G_k(x)=f(x)+\frac{1}{2}(x-x^{k-1})^TB(x-x^{k-1})-D_f(x,x^{k-1}),\end{eqnarray} which is quite similar to the \lq\lq trick\rq\rq introduced previously. However, it is not precisely the same, since the authors of \cite{BL88} do not assume that $f$ is convex, so this is not a particular case of PMAB. Unless $f$ is convex, we cannot assert that this iteration is in the SUMMA class, so we cannot be sure that the iteration reduces $\{f(x^k)\}$ to the infimal value $\beta$. This approach also relies on the extended mean value theorem, while our \lq\lq trick\rq\rq\, permits us considerable freeedom in the selection of the function $g$.

\section{Alternating Minimization (AM)}\setcounter{equation}{0}
\setcounter{theorem}{0} \setcounter{lemma}{0}
\setcounter{proposition}{0} \setcounter{corollary}{0}
\setcounter{definition}{0} \setcounter{algorithm}{0} In this section we review the basics of alternating minimization (AM) \cite{CT84}, and then show that AM, PMA and MM are equivalent. Alternating minimization plays an important role in the application of the EM algorithm \cite{DLR77} to medical image reconstruction \cite{SV82, VSK85, Byr96a}.

\subsection{The AM Method} Let $\Phi :X\times Y\rightarrow \R_+$, where $X$ and $Y$ are arbitrary nonempty sets. The objective is to find ${\hat x}\in X$ and ${\hat y}\in Y$ such that 
$$\Phi ({\hat x},{\hat y})\leq \Phi (x,y),$$ for all $x\in X$ and $y\in Y$.

The alternating minimization method \cite{CT84} is to minimize $\Phi (x,y^{k-1})$ to get $x^{k-1}$ and then to minimize $\Phi (x^{k-1},y)$ to get $y^k$. Clearly, the sequence $\{\Phi (x^{k-1},y^k)\}$ is decreasing and converges to some $\beta^*\geq -\infty$. We want $\beta^*=\Phi ({\hat x},{\hat y})$, or, at least, for $\beta^*=\beta$, where $\beta = \inf_{x,y}\Phi (x,y)$.

In AM we find $x^k$ by minimizing $\Phi (x,y^k)=\Phi (x,y(x^{k-1}))$. For each $x$ and $z$ in $X$ we define
\begin{eqnarray}d(x,z)\doteq \Phi (x,y(z))-\Phi (x,y(x)).\end{eqnarray} Clearly, $d(x,z)\geq 0$ and $d(x,x)=0$, so $d(x,z)$ is a \lq\lq distance\rq\rq . We obtain $x^k$ by minimizing
$$\Phi (x,y(x^{k-1}))=\Phi (x,y(x))+\Phi (x,y(x^{k-1}))-\Phi (x,y(x))=f(x)+d(x,x^{k-1}),$$ which shows that every AM algorithm is also a PMA. Given any AM algorithm, we define $f(x)\doteq\Phi (x,y(x))$. Then the function $g(x|z)\doteq\Phi (x,y(z))$ majorizes $f(x)$. So we see, once again, that AM, PMA and MM are equivalent methods. Now we can obtain conditions on MM algorithms sufficient for $\beta^*=\beta$ from analogous conditions expressed in the language of AM or PMA.

\subsection{The Three-Point Property} The {\it three-point property} (3PP) in \cite{CT84} is the following: for all $x\in X$ and $y\in Y$ and for all $k$ we have
\begin{eqnarray} \Phi (x,y^k)-\Phi (x^k,y^k)\geq d(x,x^k).\label{3PP}\end{eqnarray} The 3PP implies that the AM algorithm, expressed as a PMA, is in the SUMMA class and so is sufficient to have $\beta^*=\beta$.

\subsection{The Weak Three-Point Property} The 3PP is stronger than we need to get $\beta^*=\beta$; the weak 3PP implies that the AM algorithm, expressed as a PMA, is in the SUMMA2 class, and so is sufficient for $\beta^*=\beta$. The {\it weak three-point property} (w3PP) is the following: for all $x\in X$ and $y\in Y$ and for all $k$ we have
\begin{eqnarray} \Phi (x,y^k)-\Phi (x^k,y^{k+1})\geq d(x,x^k).\label{w3PP}\end{eqnarray}

\subsection{Consequences of the w3PP} From the w3PP we find that, for all $x$ and $y$,
\begin{eqnarray}d (x,x^{k-1})-d(x,x^k)\geq \Phi (x^k,y^{k+1})-\Phi (x,y(x)).\end{eqnarray} Since
$$\Phi (x^k,y^{k+1})-\Phi (x,y(x))=f(x^k)-f(x)$$ we conclude that, whenever the w3PP holds, we have
\begin{eqnarray}d(x,x^{k-1})+f(x)\geq d(x,x^k)+f(x^k),\label{summa2am}\end{eqnarray} for all $x\in X$. This means that AM with the w3PP is in the SUMMA2 class of iterative algorithms, from which it follows that $\beta^*=\beta$.

\subsection{When Do We Have $\beta^*=\beta$?} As we have noted, an AM method for which the w3PP holds is in the SUMMA2 class, so that $\beta^*=\beta$. We can formulate this in the language of MM as follows:
\begin{eqnarray}g(x|x^{k-1})-g(x|x^k)\geq f(x^k)-f(x)\end{eqnarray} for all $x$. In the language of PMA it becomes
\begin{eqnarray}d(x,x^{k-1})-d(x,x^k)\geq f(x^k)-f(x)\end{eqnarray} for all $x$.

We know that all PMAB algorithms are in the SUMMA class. Since PMA is equivalent to MM, this tells us that all MM algorithms for which $g(x|z)-f(x)$ is a Bregman distance will have $\beta^*=\beta$. As we shall see in the next section, the Auslender--Teboulle theory allows us to generalize this result.

\section{The Auslender--Teboulle Theory} In \cite{AT06} Auslender and Teboulle consider proximal minimization algorithms. They show that, if the distance $d$ has associated with it what they call \lq\lq an induced proximal distance\rq\rq  $h(x,z)$ , then $\beta^*=\beta$. It can be shown that, whenever there is an induced proximal distance, then, for any minimizer ${\hat x}$ of $f(x)$, we have
\begin{eqnarray}h({\hat x},x^k)-h({\hat x},x^{k+1})\geq f(x^k)-f({\hat x})\geq 0.\label{hh}\end{eqnarray} Consequently, the algorithm falls into the SUMMA2 class, for which $\beta^*=\beta$ is always true.

Auslender and Teboulle consider two types of distances $d$ for which there are induced proximal distances $h$: the first type are the Bregman distances, which are self-proximal in the sense that $d=h$; the second type are those having the form
$$d(x,z)=d_{\phi}(x,z)\doteq \sum_{j=1}^J z_j \phi (\frac{x_j}{z_j}),$$ for functions $\phi$ having certain properties to be discussed below. In such cases the induced proximal distance is $h(x,z)=\phi^{''}(1)KL(x,z)$, where $KL(x,z)$ is the Kullback--Leibler distance,
$$KL(x,z)=\sum_{j=1}^J x_j \log \frac{x_j}{z_j}+z_j-x_j.$$ 
Then we have
\begin{eqnarray} \phi^{''}(1)\left (KL({\hat x},x^k)-KL({\hat x},x^{k+1})\right )\geq f(x^k)-f({\hat x}).\label{phi1}\end{eqnarray} The Hellinger distance, 
$$H(x,z)=\sum_{j=1}^J (\sqrt{x_j}-\sqrt{z_j})^2,$$ fits into this framework.

The required conditions on the function $\phi (t)$ are as follows: $\phi :\R\rightarrow (-\infty ,+\infty ]$ is lower semi-continuous, proper and convex, with dom $\phi \subseteq \R_+$, and dom $\partial \phi =\R_{++}$. In addition, the function $\phi$ is $C^2$, strictly convex, and nonnegative on $\R_{++}$, with $\phi (1)=\phi '(1)=0$, and
\begin{eqnarray}\phi ''(1)\left (1-\frac{1}{t}\right )\leq \phi '(t)\leq \phi ''(1)\log (t) .\end{eqnarray}
For the Hellinger case we have $\phi (t)=(\sqrt{t}-1)^2$, so that these conditions are satisfied and we have
\begin{eqnarray}KL({\hat x},x^k)-KL({\hat x},x^{k+1})\geq 2\left (f(x^k)-f({\hat x})\right ).\label{Hell}\end{eqnarray}

We have already seen that MM algorithms for which $g(x|z)-f(x)$ is a Bregman distance have $\beta^*=\beta$. From \cite{AT06} we learn that $\beta^*=\beta$ whenever $g(x|z)-f(x)=d_{\phi}(x,z)$ for functions $\phi$ satisfying the conditions given above.

\section{AM with the Euclidean Distance}\setcounter{equation}{0}
\setcounter{theorem}{0} \setcounter{lemma}{0}
\setcounter{proposition}{0} \setcounter{corollary}{0}
\setcounter{definition}{0} \setcounter{algorithm}{0}

\subsection{Definitions} In this section we illustrate the use of AM to derive an iterative algorithm to minimize the function $f(x)=\|b-Ax\|^2$, where $A$ is an $I$ by $J$ real matrix and $b$ an $I$ by $1$ real vector. Let $R$ be the set of all $I$ by $J$ arrays $r$ with entries $r_{i,j}$ such that $\sum_{j=1}^J r_{i,j}=b_i$, for each $i$. Let $Q$ be the set of all $I$ by $J$ arrays of the form $q(x)$, where $q(x)_{i,j}=A_{i,j}x_j$. For any vectors $u$ and $v$ with the same size define
\begin{eqnarray}E(u,v)=\sum_n (u_n-v_n)^2 .\end{eqnarray}

\subsection{Pythagorean Identities}
We begin by minimizing $E(r,q(x))$ over all $r\in R$. We have the following proposition.
\begin{proposition} For all $x$ and $r$ we have
\begin{eqnarray}E(r,q(x))=E(r(x),q(x))+E(r,r(x)),\label{1stpyth}\end{eqnarray} where
\begin{eqnarray}r(x)_{i,j}=A_{i,j}x_j + \frac{1}{J}(b_i-Ax_i).\label{rx}\end{eqnarray} Therefore, $r=r(x)$ is the minimizer of $E(r,q(x))$.\end{proposition}
Now we minimize $E(r(x),q(z))$ over $z$. We have the following proposition.
\begin{proposition}For all $x$ and $z$ we have 
\begin{eqnarray}E(r(x),q(z))=E(r(x),q(Lx))+\sum_{j=1}^Jc_j(Lx_j-z_j)^2,\label{2ndpyth}\end{eqnarray} where $c_j=\sum_{i=1}^I A_{i,j}^2$ and
\begin{eqnarray}(Lx)_j=Lx_j\doteq x_j+\frac{1}{Jc_j}\sum_{i=1}^I A_{i,j}(b_i-Ax_i).\label{lx}\end{eqnarray}\end{proposition} 
We omit the proofs of these propositions, which are not deep, but involve messy calculations. Note that 
\begin{eqnarray}\|b-Ax\|^2=f(x)=JE(r(x),q(x)).\label{fxE}\end{eqnarray}

\subsection{The AM Iteration}
The iterative step of the algorithm is then
\begin{eqnarray}x^k_j=Lx^{k-1}_j=x^{k-1}_j + \frac{1}{Jc_j}\sum_{i=1}^I A_{i,j}(b_i-Ax^{k-1}_i).\label{nextx}\end{eqnarray}
Applying (\ref{1stpyth}) and (\ref{2ndpyth}) we obtain
$$f(x^{k-1})=JE(r(x^{k-1}),q(x^{k-1}))=JE(r(x^{k-1}),q(x^k))+J\sum_{j=1}^J c_j(x^k_j-x^{k-1}_j)^2$$
$$=JE(r(x^k),q(x^k))+JE(r(x^{k-1}),r(x^k))+J\sum_{j=1}^J c_j(x^k_j-x^{k-1}_j)^2$$
$$=f(x^k)+JE(r(x^{k-1}),r(x^k))+J\sum_{j=1}^J c_j(x^k_j-x^{k-1}_j)^2.$$ Therefore,
$$f(x^{k-1})-f(x^k)=JE(r(x^{k-1}),r(x^k))+J\sum_{j=1}^J c_j(x^k_j-x^{k-1}_j)^2\geq 0,$$  or
\begin{eqnarray}f(x^{k-1})-f(x^k)\geq J\sum_{j=1}^J c_j(x^k_j-x^{k-1}_j)^2\geq 0,\label{1monoE}\end{eqnarray}
from which it follows that the sequence $\{f(x^k)\}$ is decreasing and the sequence $\{\sum_{j=1}^J c_j(x^k_j-x^{k-1}_j)^2\}$ converges to zero. 

The inequality in (\ref{1monoE}) is the {\it First Monotonicity Property} for the Euclidean case. Since the sequence $\{E(b,Ax^k)\}$ is decreasing, the sequences $\{Ax^k\}$ and $\{x^k\}$ are bounded; let $x^*$ be a cluster point of the sequence $\{x^k\}$. Since the sequence $\{\sum_{j=1}^J c_j(x^k_j-x^{k-1}_j)^2\}$ converges to zero, it follows that $x^*=Lx^*$.

\subsection{Useful Lemmas} We now present several useful lemmas.

\begin{lemma} For all $x$ and $z$ we have
\begin{eqnarray}E(r(x),r(z))=\sum_{j=1}^J c_j(x_j-z_j)^2-\frac{1}{J}\sum_{i=1}^I (Ax_i-Az_i)^2.\label{err}\end{eqnarray}\end{lemma}
\begin{lemma}For all $x$ and $z$ we have
\begin{eqnarray}\frac{1}{J}\sum_{i=1}^I (Ax_i-Az_i)^2\geq \frac{1}{J^2}\sum_{j=1}^J \frac{1}{c_j}\left (\sum_{i=1}^I A_{i,j}(Ax_i-Az_i)\right )^2.\end{eqnarray}\end{lemma}
\Proof Use Cauchy's Inequality. \hfill \QED

\begin{lemma} For all $x$ and $z$ we have
\begin{eqnarray}E(r(x),r(z))\geq \sum_{j=1}^J c_j(Lx_j-Lz_j)^2.\end{eqnarray}\end{lemma}
It follows from these lemmas that this iterative algorithm is in the SUMMA2 class; for any $x$ we have
$$J\sum_{j=1}^J c_j (Lx_j-x^k)^2-J\sum_{j=1}^Jc_j (Lx_j-x^{k+1}_j)^2$$ \begin{eqnarray}\geq f(x^k)-f(x)+J\sum_{j=1}^J c_j (Lx_j-x_j)^2.\end{eqnarray}
Consequently, the sequence $f(x^k)\}$ converges to the minimum of the function $f(x)$, which must then be $f(x^*)$, and $\{x^k\}$ must converge to $x^*$.

\subsection{Characterizing the Limit} The following proposition characterizes the limit $x^*$.
\begin{proposition} The choice of ${\hat x}=x^*$ minimizes the distance $\sum_{j=1}^J c_j ({\hat x}_j-x^0_j)^2$ over all minimizers ${\hat x}$ of $f(x)=\|b-Ax\|^2$. \end{proposition}
\Proof Let ${\hat x}$ be an arbitrary minimizer of $f(x)$. Using the Pythagorean identities we find that
$$JE(r(x^k),q({\hat x}))=f({\hat x})+J\sum_{j=1}^J c_j (A{\hat x}_i-Ax^k)_i)^2-\sum_{i=1}^I (A{\hat x}_i-Ax^k_i)^2,$$ and
$$JE(r(x^k),q({\hat x}))=f(x^{k+1})+JE(r(x^k),r(x^{k+1}))+J\sum_{j=1}^J c_j ({\hat x}_j-x^{k+1}_j)^2.$$ Therefore,
$$J\sum_{j=1}^J c_j ({\hat x}_j-x^k_j)^2-J\sum_{j=1}^J c_j ({\hat x}_j-x^{k+1}_j)^2$$ $$=f(x^{k+1})-f({\hat x})+JE(r(x^k),r(x^{k+1}))+\sum_{i=1}^I (A{\hat x}_i-Ax^k_i)^2.$$
Note that the right side of the last equation depends only on $A{\hat x}$ and not directly on ${\hat x}$ itself; therefore the same is true of the left side. Now we sum both sides over the index $k$ to find that
$\sum_{j=1}^J c_j ({\hat x}_j-x^0_j)^2-\sum_{j=1}^J c_j ({\hat x}_j-x^*_j)^2$ does not depend directly on the choice of ${\hat x}$. The assertion of the proposition follows. \hfill \QED 

\subsection{SUMMA for the Euclidean Case} To get $x^k$ we minimize
$$G_k(x)=JE(r(x^{k-1}),q(x))=JE(r(x),q(x))+\left (JE(r(x^{k-1}), q(x))-JE(r(x),q(x))\right )$$ $$=f(x)+g_k(x),$$
where
$$g_k(x)=\left (JE(r(x^{k-1}), q(x))-JE(r(x),q(x))\right )=JE(r(x^{k-1}),r(x)).$$
From (\ref{err}) we have
\begin{eqnarray}g_k(x)=J\sum_{j=1}^J c_j(x^{k-1}_j-x_j)^2-\sum_{i=1}^I (Ax^{k-1}_i-Ax_i)^2.\end{eqnarray}
From
$$G_k(x)-G_k(x^k)=$$ \begin{eqnarray}JE(r(x^{k-1}),q(x))-JE(r(x^{k-1}),q(x^k))=J\sum_{j=1}^J c_j (x^k_j-x_j)^2,\end{eqnarray} we see that
$$G_k(x)-G_k(x^k)\geq g_{k+1}(x),$$ for all $x$, so that the SUMMA Inequality holds in this case. Therefore, we have
$$g_k(x)-g_{k+1}(x)\geq f(x^k)-f(x),$$ for all $x$, and so
\begin{eqnarray}g_k({\hat x})-g_{k+1}({\hat x})\geq f(x^k)-f({\hat x})\geq f(x^k)-f(x^{k+1}).\end{eqnarray} This is the {\it Second Monotonicity Property} for the Euclidean case.

\subsection{Using the Landweber Algorithm} It is of some interest to consider an alternative approach, using the Landweber (LW) algorithm. The iterative step of the LW algorithm is
\begin{eqnarray}x^k_j=x^{k-1}_j + \gamma \sum_{i=1}^I A_{i,j}(b_i-Ax^{k-1}_i),\end{eqnarray} where $0<\gamma <\frac{2}{\rho (A^TA)}$.
We define $\beta_j=\frac{1}{Jc_j}$, $B_{i,j}=\sqrt{\beta_j}A_{i,j}$, and $z_j=x_j/\sqrt{\beta_j}$. Then $Bz=Ax$. The LW algorithm, applied to $Bz=b$ and with $\gamma =1$,
is
\begin{eqnarray}z^k=z^{k-1}+B^T(b-Bz^{k-1}).\end{eqnarray} Since the trace of $B^TB$ is one, the choice of $\gamma =1$ is allowed. It is known that the LW algorithm converges to the minimizer of $\|b-Bz\|$ for which $\|z-z^0\|$ is minimized. Converting back to the original $x^k$, we find that we get the same iterative sequence that we got using the AM method. Moreover, we find once again that the sequence $\{x^k\}$ converges to the minimizer $x^*$ of $f(x)$ for which the distance $\sum_{j=1}^J c_j ({\hat x}_j-x^0_j)^2$ is minimized over all minimizers ${\hat x}$ of $f(x)$.

The Landweber algorithm applied to the original problem of minimizing $f(x)=\|Ax-b\|^2$ has the iterative step
\begin{eqnarray}x^k=x^{k-1}-\gamma A^T(Ax^{k-1}-b),\end{eqnarray} where $0<\gamma <\frac{2}{\rho (A^TA)}$. The sequence $\{x^k\}$ converges to the minimizer $x^*$ of $f(x)$ that minimizes $\|{\hat x}-x^0\|$ over all minimizers ${\hat x}$ of $f(x)$. 

\section{The SMART}\setcounter{equation}{0}
\setcounter{theorem}{0} \setcounter{lemma}{0}
\setcounter{proposition}{0} \setcounter{corollary}{0}
\setcounter{definition}{0} \setcounter{algorithm}{0} 

In this section we discuss the {\it simultaneous multiplicative algebraic reconstruction technique} (SMART) \cite{DR72, Sch72, CS87,Byr93, Byr95, Byr96a}. A key step in the proof of convergence is showing that the SMART is in the SUMMA class.

\subsection{The Kullback--Leibler or Cross-Entropy Distance} The Kullback--Leibler distance is quite useful in the discussions that follow. For positive numbers $s$ and $t$, the Kullback--Leibler distance from $s$ to $t$ is
\begin{eqnarray}KL(s,t)=s\log \frac{s}{t} +t-s.\label{klst}\end{eqnarray} Since, for $x>0$ we have
$$x-1-\log x\geq 0$$ and equal to zero if and only if $x=1$, it follows that
$$KL(s,t)\geq 0,$$ and $KL(s,s)=0$. We use limits to define $KL(0,t)=t$ and $KL(s,0)=+\infty$. Now we extend the KL distance to nonnegative vectors component-wise. The following lemma is easy to prove.
\begin{lemma}\label{klsumlemma} For any nonnegative vectors $x$ and $z$, with $z_+=\sum_{j=1}^J z_j>0$, we have 
\begin{eqnarray} KL(x,z) = KL(x_+,z_+)+ KL(x,\frac{x_+}{z_+}z).\label{klsum}\end{eqnarray} \end{lemma} We can extend the KL distance in the obvious way to infinite sequences with nonnegative terms, as well as to nonnegative functions of continuous variables.
 
\subsection{The Problem to be Solved} We assume that $y$ is a positive vector in $\R^I$, $P$ an $I$ by $J$ matrix with nonnegative entries $P_{i,j}$, $s_j=\sum_{i=1}^IP_{i,j}>0$, and we want to find a nonnegative solution or approximate solution $x$ for the linear system of equations $y=Px$. The SMART will minimize $KL(Px,y)$, over $x\geq 0$. For notational simplicity we shall assume that the system has been normalized so that $s_j=1$ for each $j$.

\subsection{The SMART Iteration} The SMART algorithm \cite{DR72, Sch72, CS87, Byr93, Byr96a} minimizes the function
$f(x)=KL(Px,y)$, over nonnegative vectors $x$. Having found the vector $x^{k-1}$, the next vector in the SMART
sequence is $x^k$, with entries given by
\begin{eqnarray}x^k_j=x^{k-1}_j\exp \Big (\sum_{i=1}^I P_{ij}\log
(y_i/(Px^{k-1})_i)\Big ).\label{smart2}\end{eqnarray}
The iterative step of the SMART can be decsribed as $x^k=Sx^{k-1}$, where $S$ is the operator defined by
\begin{eqnarray}(Sx)_j\doteq x_j\exp \Big (\sum_{i=1}^I P_{ij}\log
(y_i/(Px)_i)\Big ).\label{smartT}\end{eqnarray}
In our proof of convergence of the SMART we will show that any cluster point $x^*$ of the SMART sequence $\{x^k\}$ is a fixed point of the operator $S$. To avoid pathological cases in which $Px^*_i=0$ for some index $i$, we can assume, at the outset, that all the entries of $P$ are positive. This is wise, in any case, since the model of $y=Px$ is unlikely to be exactly accurate in applications. 

\subsection{The SMART as AM} In \cite{Byr93} the SMART was
derived using the following alternating minimization (AM) approach. 

For each $x$, let $r(x)$ and $q(x)$ be the $I$ by $J$
arrays with entries
\begin{eqnarray}r(x)_{ij}=x_jP_{ij}y_i/(Px)_i ,\end{eqnarray}
and
\begin{eqnarray}q(x)_{ij}=x_jP_{ij}.\end{eqnarray}

In the iterative step of the SMART we get $x^k$ by minimizing the
function

\begin{eqnarray}G_k(x)=KL(q(x),r(x^{k-1}))=\sum_{i=1}^I\sum_{j=1}^J KL(q(x)_{ij},r(x^{k-1})_{ij})\label{GkS}\end{eqnarray}
over $x\geq 0$. Note that $f(x)=KL(Px,y)=KL(q(x),r(x))$.
We have the following helpful {\it Pythagorean identities}:
\begin{eqnarray}KL(q({ x}),r({ z}))=KL(q({ x}),r({ x}))+KL({ x},{ z})-KL(P{ x},P{ z});\label{py3}\end{eqnarray} and
\begin{eqnarray}KL(q({ x}),r({ z}))=KL(q({ Sz}),r({ z}))+KL({ x},{ Sz}).\label{py4}\end{eqnarray}
Note that it follows from Equation (\ref{klsum}) that $KL({ x},{ z})-KL(P{ x},P{
z})\geq 0$.

From the Pythagorean identities we find that $x^k$ is obtained by minimizing
\begin{eqnarray}G_k(x)=KL(Px,y)+KL(x,x^{k-1})-KL(Px,Px^{k-1}),\end{eqnarray} so that SMART is an AF algorithm and 
\begin{eqnarray}g_k(x)=KL(x,x^{k-1})-KL(Px,Px^{k-1}).\label{gkS}\end{eqnarray} Consequently, the sequence $\{f(x^k)\}$ is decreasing and the sequences $\{Px^k\}$ and $\{x^k\}$ are bounded. From
$$G_k(x)-G_k(x^k)=KL(x,x^k)\geq KL(x,x^k)-KL(Px,Px^k)=g_{k+1}(x)$$ we conclude that the SMART is in the SUMMA class.
It follows from our discussion of the SUMMA Inequality that, for all $x\geq 0$,
\begin{eqnarray}g_k(x)+f(x)\geq g_{k+1}(x)+f(x^k).\label{SMARTSUMMA1}\end{eqnarray} Since
$$\sum_{j=1}^J x_j^k \leq \sum_{i=1}^I y_i,$$ we see once again that the sequence $\{x^k\}$ is bounded and therefore has a cluster point, $x^*$, with $f(x^k)\geq f(x^*)$ for all $k$ and $Sx^*=x^*$. 

\subsection{MM in SMART} At each step of the SMART we minimize the function $KL(q(x),r(x^{k-1}))$ to get $x^k$. From
\begin{eqnarray}KL(q(x),r(z))=KL(q(x),r(x))+KL(x,z)-KL(Px,Pz)\geq KL(Px,y)\end{eqnarray} we see that the function $KL(q(x),r(z))=g(x|z)$ is a majorizing function for the function $f(x)=KL(Px,y)$. 

\subsubsection{The First Monotonicity Property for SMART} Using the Pythagorean identities we have
\begin{eqnarray}KL(Px^k,y)-KL(Px^{k+1},y)\geq KL(x^k,x^{k+1}).\label{1ms}\end{eqnarray}

\subsubsection{The Second Monotonicity Property for SMART}  Let ${\hat x}$ be any minimizer of $KL(Px,y)$. We then have
\begin{eqnarray}KL({\hat x},x^k)-KL({\hat x},x^{k+1}) =KL(Px^{k+1},y)-KL(P{\hat x},y)+ \cr KL(P{\hat x},Px^k)+KL(x^{k+1},x^k)-KL(Px^{k+1},Px^k)\geq 0.\label{2ms}\end{eqnarray}

In fact, there is a somewhat more general version of (\ref{2ms}), that tells us that, since $Sx^*=x^*$ and $f(x^k)\geq f(x^*)$, we can replace ${\hat x}$ with $x^*$ in (\ref{2ms}), to get 
\begin{eqnarray}KL(x^*,x^k)-KL(x^*,x^{k+1}) =KL(Px^{k+1},y)-KL(Px^*,y)+ \cr KL(Px^*,Px^k)+KL(x^{k+1},x^k)-KL(Px^{k+1},Px^k)\geq 0.\label{2ms2}\end{eqnarray}

From (\ref{2ms2}) it follows that the sequence $\{f(x^k)\}$ converges to $f(x^*)$. Since the SMART is in SUMMA, we know that $f(x^*)$ must be the minimum of $f(x)$. Since a subsequence of $\{KL(x^*,x^k)\}$ converges to zero, it follows that $\{x^k\}$ converges to $x^*$.

\subsection{Characterizing the Limit of SMART} Let ${\hat x}$ be any minimizer of $KL(Px,y)$.
From Equation (\ref{2ms}) we see that the difference $KL({\hat x},x^k)-KL({\hat x},x^{k+1})$ depends only on $P{\hat x}$, and not on ${\hat x}$ itself. Summing over the index $k$ on both sides and \lq\lq telescoping\rq\rq\,, we find that the difference $KL({\hat x},x^0)-KL({\hat x},x^*)$ also depends only on $P{\hat x}$, and not on ${\hat x}$ itself. It follows that ${\hat x}=x^*$ is the minimizer of $f(x)$ for which $KL({\hat x},x^0)$ is minimized. If $y=Px$ has nonnegative solutions, and the entries of $x^0$ are all equal to one, then $x^*$ maximizes the Shannon entropy over all nonnegative solutions of $y=Px$.

The following theorem summarizes the situation with regard to the
SMART \cite{Byr93, Byr95, Byr96a}.

\begin{theorem}In the consistent case, in which the system $y=Px$ has nonnegative solutions, the sequence of iterates of SMART converges to
the unique nonnegative solution of $y=Px$ for which the distance
$KL(x,x^0)$ is minimized. In the inconsistent
case it converges to the unique nonnegative minimizer of the
distance $KL(Px,y)$ for which $KL(x,x^0)$ is
minimized. In the inconsistent case, if $P$ and every matrix derived from $P$ by deleting
columns has full rank then there is a unique nonnegative minimizer
of $KL(Px,y)$ and at most $I-1$ of its entries are nonzero.
\end{theorem}

\section{The EMML Algorithm}\setcounter{equation}{0}
\setcounter{theorem}{0} \setcounter{lemma}{0}
\setcounter{proposition}{0} \setcounter{corollary}{0}
\setcounter{definition}{0} \setcounter{algorithm}{0} 

In this section we discuss the EMML algorithm \cite{SV82, VSK85,Byr93, Byr95, Byr96a}. A key step in the proof of convergence is showing that the EMML algorithm is in the SUMMA2 class.

\subsection{The EMML Iteration}  Once again, we want to find a nonnegative solution or approximate solution $x$ for the linear system of equations $y=Px$. The EMML algorithm will minimize $KL(y,Px)$.

The EMML algorithm\index{EMML algorithm} minimizes the function
$f(x)=KL(y,Px)$, over nonnegative vectors $x$. Having found the
vector $x^{k-1}$, the next vector in the EMML sequence is $x^k$,
with entries given by
\begin{eqnarray}x^k_j=x^{k-1}_j\Big (\sum_{i=1}^I P_{ij}
(y_i/(Px^{k-1})_i)\Big ).\label{emmlb}\end{eqnarray}
The iterative step of the EMML algorithm can be described as $x^k=Mx^{k-1}$, where $M$ is the operator defined by
\begin{eqnarray}(Mx)_j\doteq x_j\Big (\sum_{i=1}^I P_{ij}
(y_i/(Px)_i)\Big ).\label{emmlS}\end{eqnarray}
As we shall see, the EMML algorithm forces the sequence $\{KL(y,Px^k)\}$ to be decreasing. It follows that $(Px^*)_i>0$, for any cluster point $x^*$ and for all $i$.

\subsection{The EMML  as AM} Now we want to minimize $f(x)=KL(y,Px)$.
We have the following helpful {\it Pythagorean identities}:
\begin{eqnarray}KL(r({ x}),q({ z}))=KL(r({ z}),q({ z}))+KL(r({ x}),r({ z}));\label{py1}\end{eqnarray} and 
\begin{eqnarray}KL(r({ x}),q({ z}))=KL(r({ x}),q({ Mx}))+KL({ Mx},{ z}).\label{py2}\end{eqnarray}
Using these Pythagorean identities we see that, for $\{{ x}^k\}$ given by Equation (\ref{emmlb}), the sequence $\{KL({y},P{ x}^k)\}$ is decreasing and the sequences $\{KL({ x}^{k+1},{ x}^k)\}$ and $\{KL(r(x^k),r(x^{k+1}))\}$ converge to zero. It follows that the EMML sequence $\{{ x}^k\}$ is bounded. 
In fact, we have
$$\sum_{j=1}^J x^k_j = \sum_{i=1}^I y_i.$$
Using (\ref{klsum}) we obtain the following useful inequality:
\begin{eqnarray}KL(r(x),r(z))\geq KL(Mx,Mz).\label{rxMx}\end{eqnarray} From
$$KL(r(x),q(x^k))=KL(r(x^k),q(x^k))+KL(r(x),r(x^k))\geq f(x^k)+KL(Mx,x^{k+1}),$$ and
$$KL(r(x),q(x^k))=KL(r(x),q(Mx))+KL(Mx,x^k)=f(x)-KL(Mx,x)+KL(Mx,x^k)$$ we have
\begin{eqnarray}KL(Mx,x^k)-KL(Mx,x^{k+1})\geq f(x^k)-f(x)+KL(Mx,x).\label{s24}\end{eqnarray} Note that we have used (\ref{rxMx}) here. Therefore, the EMML is in the SUMMA2 class. With $x^*$ a cluster point, we have
\begin{eqnarray}KL(Mx^*,x^k)-KL(Mx^*,x^{k+1})\geq f(x^k)-f(x^*)\geq 0.\label{s24a}\end{eqnarray}
Therefore, the sequence $\{KL(Mx^*,x^k)\}$ is decreasing, and the sequence $\{f(x^k)\}$ converges to $f(x^*)$. Since the EMML is in the SUMMA2 class, we know that $f(x^*)$ is the minimum value of $f(x)$ and $Mx^*=x^*$.

The following theorem summarizes the situation with regard to the
EMML algorithm \cite{Byr93, Byr95, Byr96a}.

\begin{theorem}In the consistent case, in which the system $y=Px$ has nonnegative solutions, the sequence of EMML iterates converges to
a nonnegative solution of $y=Px$. In the inconsistent
case it converges to a nonnegative minimizer of the
distance $KL(y,Px)$. In the inconsistent case, if $P$ and every matrix derived from $P$ by deleting
columns has full rank then there is a unique nonnegative minimizer
of $KL(y,Px)$ and at most $I-1$ of its entries are nonzero.
\end{theorem} In contrast with the SMART, we have been unable to characterize the limit in terms of the starting vector $x^0$.

\subsection{MM in EMML} At each step of the EMML algorithm we minimize $KL(r(x^{k-1}),q(x))$ to get $x^k$. From
\begin{eqnarray}KL(r(z),q(x))=KL(r(x),q(x))+KL(r(z),r(x))\end{eqnarray}
we see that the function
\begin{eqnarray}KL(r(z),q(x))=g(x|z)\end{eqnarray} is a majorizing function for $f(x)=KL(y,Px)$. 
\subsection{The First Monotonicity Property for EMML} From the Pythagorean identities we have
\begin{eqnarray}KL(y,Px^k)-KL(y,Px^{k+1})=KL(r(x^k),r(x^{k+1}))+KL(x^{k+1},x^k),\end{eqnarray} so that
\begin{eqnarray}KL(y,Px^k)-KL(y,Px^{k+1})\geq KL(x^{k+1},x^k).\label{firstmono}\end{eqnarray} The inequality in (\ref{firstmono}) is called the {\it First Monotonicity Property} in \cite{EL98}. 
\subsection{The Second Monotonicity Property for EMML} Let ${\hat x}$ be a minimizer of $f(x)=KL(y,Px)$. Inserting $x={\hat x}$ into Equation (\ref{s24}), we obtain
\begin{eqnarray}KL({\hat x},x^k)-KL({\hat x},x^{k+1})\geq KL(y,Px^k)-KL(y,Px^{k+1}).\label{secondmono}\end{eqnarray} The inequality in (\ref{secondmono}) is called the {\it Second Monotonicity Property} in \cite{EL98}.

\section{The Hellinger Distance}\setcounter{equation}{0}
\setcounter{theorem}{0} \setcounter{lemma}{0}
\setcounter{proposition}{0} \setcounter{corollary}{0}
\setcounter{definition}{0} \setcounter{algorithm}{0}In \cite{EL98} the authors consider extending the results concerning the KL distance to the Hellinger distance. In particular, they explore the use of AM and MM. 

\subsection{The Definition of $H(s,t)$} For $s>0$ and $t>0$ the Hellinger distance from $s$ to $t$ is
\begin{eqnarray}H(s,t)= (\sqrt{s}-\sqrt{t})^2.\label{H}\end{eqnarray} As in the case of the KL distance, we can extend $H$ to nonnegative vectors component-wise. In this section we consider the problem of minimizing $H(y,Px)$ given by
\begin{eqnarray}H(y,Px)=\sum_{i=1}^I (\sqrt{y_i}-\sqrt{(Px)_i})^2.\label{H2}\end{eqnarray}
As in the KL case, we assume that $s_j=\sum_{i=1}^I P_{i,j}=1$ for each $j$.

\subsection{An Alternating-Minimization Approach} 

In (4.2) of \cite{EL98} the authors present a majorizing function to be used to generate the iterative sequence. 
We can show that their majorizing function is $g(x|z)=H(r(z),q(x))$, with the same notation as in the KL case. 
The following proposition is essentially what appears in \cite{EL98}. The proof here is simpler than in \cite{EL98}.
\begin{proposition} For all $x\geq 0$ and $z\geq 0$ we have
\begin{eqnarray}\sum_{j=1}^J P_{i,j}\sqrt{x_jz_j}\leq \sqrt{(Px)_i (Pz)_i}.\label{ellin}\end{eqnarray}
\end{proposition}
\Proof We have
$$\sum_{j=1}^J P_{i,j}\sqrt{x_j}\sqrt{z_j}=\sum_{j=1}^J \sqrt{P_{i,j}x_j}\sqrt{P_{i,j}z_j}$$
$$\leq \sqrt{\sum_{j=1}^J P_{i,j}x_j}\sqrt{\sum_{j=1}^J P_{i,j}z_j}=\sqrt{(Px)_i(Pz)_i},$$ by the Cauchy Inequality. \hfill \QED 

\begin{corollary} The function $g(x|z)=H(r(z),q(x))$ majorizes $f(x)=H(y,Px)$. \label{Hmajor}\end{corollary}
\Proof We only need to show that
$$\sum_{i=1}^I \sum_{j=1}^J \sqrt{r(z)_{i,j}q(x)_{i,j}}\geq \sum_{i=1}^I \sqrt{(Px)_i y_i}.$$
We have
$$\sqrt{r(z)_{i,j}q(x)_{i,j}}=P_{i,j}\sqrt{z_j x_j}\sqrt{\frac{y_i}{(Pz)_i}},$$
from which it follows that
$$\sum_{i=1}^I \sum_{j=1}^J \sqrt{r(z)_{i,j}q(x)_{i,j}}=\sum_{i=1}^I \sum_{j=1}^J P_{i,j}\sqrt{z_j x_j}\sqrt{\frac{y_i}{(Pz)_i}}$$ 
$$=\sum_{i=1}^I \left (\sum_{j=1}^J P_{i,j}\sqrt{z_j x_j}\right )\sqrt{\frac{y_i}{(Pz)_i}}$$
$$\leq  \sum_{i=1}^I \left (\sqrt{(Pz)_i (Px)_i}\right )\sqrt{\frac{y_i}{(Pz)_i}}=\sum_{i=1}^I \sqrt{(Px)_i y_i}.$$ \hfill \QED

Note that Corollary \ref{Hmajor} can also be obtained by using Lagrange multipliers to minimize $H(r,q(x))$ over all $r=\{r_{i,j}\}$ with $\sum_{j=1}^J r_{i,j}=y_i$, for all $i$.

\begin{corollary} For all $x\geq 0$ and $z\geq 0$ we have
\begin{eqnarray}\sum_{i=1}^I \sqrt{(Px)_i (Pz)_i} \geq \sum_{j=1}^J \sqrt{x_j z_j}.\label{cor2}\end{eqnarray} \end{corollary}
\Proof From
$$\sum_{j=1}^J P_{i,j}\sqrt{x_jz_j}\leq \sqrt{(Px)_i (Pz)_i}$$ we have
$$\sum_{i=1}^I \sum_{j=1}^J P_{i,j}\sqrt{x_jz_j}\leq \sum_{i=1}^I \sqrt{(Px)_i (Pz)_i},$$ so that
$$\sum_{j=1}^J \left (\sum_{i=1}^I P_{i,j}\right )\sqrt{x_jz_j}\leq \sum_{i=1}^I \sqrt{(Px)_i (Pz)_i}.$$ \hfill \QED 

The iterative step of the algorithm is derived by minimizing $H(r(x^{k-1}),q(x))$ to get $x^k$, with
\begin{eqnarray}x^k_j=x^{k-1}_j\left (\sum_{i=1}^I P_{i,j}\frac{\sqrt{y_i}}{\sqrt{(Px^{k-1})_i}}\right )^2.\end{eqnarray}
We can write $x^k=Tx^{k-1}$, where $T$ is the operator
\begin{eqnarray}Tx_j=x_j\left (\sum_{i=1}^I P_{i,j}\frac{\sqrt{y_i}}{\sqrt{(Px)_i}}\right )^2.\label{T}\end{eqnarray}
Since $g(x|z)$ majorizes $f(x)$, it follows easily that the sequence $\{f(x^k)\}$ is decreasing, so that the sequences $\{Px^k\}$ and $\{x^k\}$ are bounded. 

In the EMML case we saw that 
$$\sum_{i=1}^I y_i=\sum_{j=1}^J Mx_j,$$ while for SMART we have
$$\sum_{i=1}^I y_i \geq \sum_{j=1}^J Sx_j.$$  In the Hellinger case we shall see that
$$\sum_{i=1}^I y_i \geq \sum_{j=1}^J Tx_j.$$

In the EMML case we have the Pythagorean identity
\begin{eqnarray}KL(r(x),q(z))=KL(r(x),q(Mx))+KL(Mx,z),\end{eqnarray} while in the Hellinger case we have the analogous Pythagorean identity
\begin{eqnarray}H(r(x),q(z))=H(r(x),q(Tx))+H(Tx,z),\end{eqnarray} so that
\begin{eqnarray}H(r(x),q(x))=H(r(x),q(Tx))+H(Tx,x).\end{eqnarray}
We note that, unlike the KL distance, the Hellinger distance is symmetric; we have
\begin{eqnarray}H(x,z)=H(z,x).\end{eqnarray}

\begin{lemma} For every $x\geq 0$ we have 
\begin{eqnarray}H(r(x),q(Tx))=\sum_{i=1}^I y_i-\sum_{j=1}^J (Tx)_j \geq 0,\end{eqnarray} so that the set $\{Tx| x\geq 0\}$ is bounded. \end{lemma} Since minimizing $f(x)=H(r(x),q(x))$ is equivalent to minimizing $H(r(x),q(Tx))$, it follows that minimizing $f(x)$ is equivalent to maximizing $\sum_{j=1}^J Tx_j$.
In the EMML case we have
\begin{eqnarray}f(x)=KL(y,Px)=KL(r(x),q(x)),\end{eqnarray} while in the Hellinger case we have the analogous result
\begin{eqnarray}f(x)=H(y,Px)=H(r(x),q(x)).\end{eqnarray}
In the EMML case we use the other Pythagorean identity
\begin{eqnarray}KL(r(z),q(x))=KL(r(x),q(x))+KL(r(z),r(x))\label{emmaj}\end{eqnarray} to show that
\begin{eqnarray}KL(r(z),q(x))=g(x|z)\end{eqnarray} is a majorizing function for $f(x)=KL(y,Px)$. In the Hellinger case we have shown that
\begin{eqnarray}H(r(z),q(x))\geq H(r(x),q(x)),\end{eqnarray} for all $x\geq 0$. It would be nice if we had an analogue of Equation (\ref{emmaj}) for the Hellinger case. Said another way, can we find a simple expression for
$$H(r(z),q(x))-H(r(x),q(x))?$$ By analogy with the EMML case, we might expect to have
\begin{eqnarray}H(r(z),q(x))-H(r(x),q(x))=H(r(z),r(x)).\label{Hrr}\end{eqnarray} Actually, we don't need this much; it would be enough to prove that Equation (\ref{Hrr}) holds for $x=Tz$.

It is worth noting here that perhaps we should consider analogies, not just with the EMML, but with the SMART also. The Hellinger distance is symmetric, so that $H(y,Px)=H(Px,y)$, whereas $KL(y,Px)$ and $KL(Px,y)$ are not the same. In the SMART case we have the inequality
\begin{eqnarray}KL(x,z)\geq KL(Px,Pz).\end{eqnarray} This holds as well for the Hellinger distance.
\begin{lemma} For all $x\geq 0$ and $z\geq 0$ we have 
\begin{eqnarray}H(x,z)\geq H(Px,Pz).\label{HxPx}\end{eqnarray}\end{lemma}
\Proof We have
$$H(Px,Pz)=\sum_{i=1}^I \left ((Px)_i +(Pz)_i-2 \sqrt{(Px)_i (Pz)_i}\right )$$
$$=\sum_{j=1}^J\left (x_j + z_j \right ) -2\sum_{i=1}^I \sqrt{(Px)_i (Pz)_i}$$
$$\leq \sum_{j=1}^J\left (x_j + z_j \right ) -2\sum_{j=1}^J \sqrt{x_j z_j},$$
by (\ref{cor2}). \hfill \QED

\subsection{Convergence}  From the discussion above we have
\begin{eqnarray}H(y,Px^k)-H(y,Px^{k+1})\geq H(x^k,x^{k+1}),\end{eqnarray} so that the sequence $\{H(y,Px^k)\}$ is decreasing and the sequence $\{H(x^k,x^{k+1})\}$ converges to zero. Since the sequence $\{x^k\}$ is bounded, it has a cluster point, call it ${\hat x}$ which must then be a fixed point of $T$. The sequence $\{H(y,Px^k)\}$ then converges to $H(y,P{\hat x})$. In \cite{EL98} it was shown that, if ${\hat x}$ minimizes $f(x)$, then
\begin{eqnarray}KL({\hat x},x^k)-KL({\hat x},x^{k+1})\geq 2\left (H(y,Px^{k+1})-H(y,P{\hat x})\right ).\label{KLH}\end{eqnarray} It follows that the sequence $\{x^k\}$ converges to ${\hat x}$.

\section{Pearson's $\phi^2$ Distance}\setcounter{equation}{0}
\setcounter{theorem}{0} \setcounter{lemma}{0}
\setcounter{proposition}{0} \setcounter{corollary}{0}
\setcounter{definition}{0} \setcounter{algorithm}{0} In \cite{EL98} the authors consider extending the results concerning the KL and Hellinger distances to the $\phi^2$-distance of Pearson . In particular, they explore the use of AM and MM.

\subsection{The Definition of $\phi^2(s,t)$} For $s>0$ and $t>0$ Pearson's $\phi^2$ distance from $s$ to $t$ is
\begin{eqnarray}\phi^2(s,t)=\frac{(s-t)^2}{t}\label{P}\end{eqnarray} As in the cases of the KL and H distances, we can extend $\phi^2$ to nonnegative vectors component-wise. Note that $\phi^2(s,t)$ is not symmetric. In this section we consider the problem of minimizing $\phi^2(y,Px)$ given by
\begin{eqnarray}\phi^2(y,Px)=\sum_{i=1}^I \frac{(y_i-(Px)_i)^2}{(Px)_i}.\label{P2}\end{eqnarray}
As in the previous cases, we assume that $s_j=1$ for each $j$.

\subsection{An Alternating-Minimization Approach} 
In (5.4) of \cite{EL98} the authors present a majorizing function to be used to generate the iterative sequence. We can show that their majorizing function is $g(x|z)=\phi^2(r(z),q(x))$, with the same notation as in the KL and H  cases. The following proposition is essentially what appears in \cite{EL98}; the proof given here is simpler, however.
\begin{proposition} For all $x>0$ and $z>0$ we have
\begin{eqnarray}\sum_{j=1}^J P_{i,j}\frac{x_j^2}{z_j} \geq  \frac{(Px)_i^2}{(Pz)_i}.\label{ellin2}\end{eqnarray}
\end{proposition}
\Proof We have
$$(Px)_i=\sum_{j=1}^J P_{i,j}x_j=\sum_{j=1}^J \sqrt{P_{i,j}z_j}\sqrt{P_{i,j}\frac{x_j^2}{z_j}}$$
$$\leq \sqrt{\sum_{j=1}^J P_{i,j}z_j}\sqrt{\sum_{j=1}^J P_{i,j}\frac{x_j^2}{z_j}},$$ so that
$$(Px)_i^2\leq (Pz)_i\sum_{j=1}^J P_{i,j}\frac{x_j^2}{z_j}.$$
\hfill \QED

\begin{corollary} For all $x>0$ and $z>0$ we have
\begin{eqnarray}\sum_{i=1}^I \frac{(Px)_i^2}{(Pz)_i} \leq \sum_{j=1}^J \frac{x_j^2}{z_j} .\end{eqnarray}
\end{corollary}

\begin{corollary} The function $g(x|z)=\phi^2(r(z),q(x))$ majorizes $\phi^2(y,Px)$. \label{Pmajor}\end{corollary} 
Note that Corollary \ref{Pmajor} can also be obtained by using Lagrange multipliers to minimize $\phi^2(r,q(x))$ over all $r=\{r_{i,j}\}$ with $\sum_{j=1}^J r_{i,j}=y_i$, for all $i$.
\begin{corollary} For each $x>0$ and $z>0$ we have
\begin{eqnarray} \phi^2(x,z)\geq \phi^2(Px,Pz).\end{eqnarray}
\end{corollary}

The iterative step of the algorithm is derived by minimizing $\phi^2(r(x^{k-1}),q(x))$ to get $x^k$ given by
\begin{eqnarray}x_j^k=x_j^{k-1}\sqrt{\sum_{i=1}^IP_{i,j}\left (\frac{y_i}{(Px^{k-1})_i}\right )^2}.\end{eqnarray}
With $R$ the operator defined by
\begin{eqnarray}(Rx)_j=Rx_j\doteq x_j\sqrt{\sum_{i=1}^IP_{i,j}\left (\frac{y_i}{(Px)_i}\right )^2},\end{eqnarray} we can write $x^k=Rx^{k-1}$.
An easy calculation shows that $\phi^2(Rz,x)=\phi^2(q(Rz),q(x))$ and
\begin{eqnarray}\phi^2(r(z),q(x))=\phi^2(r(z),q(Rz))+\phi^2(Rz,x).\end{eqnarray} Since $g(x|z)$ majorizes $f(x)$ it follows that the sequence $\{f(x^k)\}$ is decreasing, so that the sequences $\{Px^k\}$ and $\{x^k\}$ are bounded. We also have
$$\phi^2(r(x),q(Rx))=\sum_{j=1}^J (Rx)_j-\sum_{i=1}^I y_i \geq 0.$$

\section{Just a Coincidence?}\setcounter{equation}{0}
\setcounter{theorem}{0} \setcounter{lemma}{0}
\setcounter{proposition}{0} \setcounter{corollary}{0}
\setcounter{definition}{0} \setcounter{algorithm}{0}

As we have seen, the KL distance appears, apparently uninvited, in (\ref{KLH}). In \cite{AT06} a similar thing happens, as (\ref{Hell}) shows, prompting us to ask if this is just a coincidence, or if something deeper is going on here.

In proximal minimization algorithms (PMA) we obtain an iterative method for minimizing a function $f(x)$ by minimizing
$$f(x)+d(x,x^{k-1})$$ to get the next iterate $x^k$. Here $d(x,z)\geq 0$ and $d(x,x)=0$, for all $x$ and $z$. It follows easily that the sequence $\{f(x^k)\}$ is decreasing to a limit $\beta^*\geq -\infty$. We have discussed what additional restrictions should be placed on the distance $d$ to guarantee that 
$$\beta^*=\beta \doteq \inf_x \{f(x)\}.$$ 
For the Hellinger distance we have $H(x,z)=d_{\phi}(x,z)$, for $\phi (t)=(\sqrt{t}-1)^2$, so that, according to \cite{AT06},
\begin{eqnarray}KL({\hat x},x^k)-KL({\hat x},x^{k+1})\geq 2\left (f(x^k)-f({\hat x})\right ).\label{Hell}\end{eqnarray} This looks a lot like (\ref{KLH}).

Of course, the problems are not quite the same; in \cite{AT06} they are trying to minimize some unrelated function $f(x)$, using the Hellinger distance in the PMA framework, while we are trying to minimize $H(y,Px)=H(Px,y)$ using alternating minimization. However, the resemblance between (\ref{KLH}) and (\ref{Hell}) must be more than a coincidence, mustn't it?

\section{Acknowledgments} The authors thank Professor Paul Eggermont of the University of Delaware and Professor Hung Phan of the University of Massachusetts Lowell for helpful discussions, and Professor Eggermont for making available the preprint \cite{EL98}.

\makeatletter
\renewcommand\@biblabel[1] {#1.}
\makeatother

\end{document}